\documentclass[paper=a4, 12pt]{article} 

\usepackage[utf8]{inputenc} 
\usepackage{geometry}
\usepackage{amsmath,amsfonts,amsthm,amssymb,mathtools} 
\usepackage{enumerate}
\usepackage{stmaryrd}
\usepackage[english]{babel}
\usepackage{relsize}
\usepackage{color}
\usepackage{fancyhdr}
 \usepackage{multicol}
 \usepackage{tikz}
 \usepackage[all]{xy}
 \usepackage[noadjust]{cite}
 \usepackage[hyperfootnotes=false]{hyperref}

\usepackage{url}
\SelectTips{cm}{12}

\numberwithin{figure}{section} 
\numberwithin{table}{section} 

\theoremstyle{theorem}
\newtheorem{theorem}{Theorem}[section]
\newtheorem{proposition}[theorem]{Proposition}

\newtheorem{lemma}[theorem]{Lemma}

\theoremstyle{definition}

\newtheorem{remark}[theorem]{Remark}

\newcommand{\ulpullback}[1][ul]{\save*!/#1-4ex/#1:(-1,1)@^{|-}\restore}

\newcommand{\drpullback}[1][dr]{\save*!/#1-4ex/#1:(-1,1)@^{|-}\restore} 
\newcommand{\dpullbackS}[1][d]{\save*!/#1-3ex/#1:(-1,1)@^{|-}\restore}
\newcommand{\dpullback}[1][d]{\save*!/#1-4ex/#1:(-1,1)@^{|-}\restore}

\DeclareMathAlphabet{\mathbbe}{U}{bbold}{m}{n}
\newcommand{\simplexcategory}{\mathbbe{\Delta}} 

\newcommand{\pow}[1]{\llbracket #1\rrbracket}

\newcommand{\infgrpd}{\mathbf{Grpd}}
\newcommand{\sur}{\mathbf{S}}

\newcommand{\grpdSprod}{\odot}

\newcommand{\arxiv}[1]{\href{http://arxiv.org/pdf/#1}{arXiv:#1}}

\newcommand{\Q}{\mathbb Q}
\newcommand{\name}[1]{\ulcorner #1\urcorner}

\newcommand{\C}{\mathcal C}
\newcommand{\B}{\mathbb B}
\newcommand{\Sg}{\mathbb S}
\newcommand{\T}{\mathbb T}
\newcommand{\z}{\mathbf z}
\newcommand{\x}{\mathbf x}
\renewcommand{\l}{\lambda}

\DeclareMathOperator{\autiv}{aut}
\DeclareMathOperator{\aut}{Aut}
\DeclareMathOperator{\rep}{R}
\DeclareMathOperator{\iso}{Iso}
\DeclareMathOperator{\tw}{Tw^+}
\DeclareMathOperator{\trans}{Trans}
\DeclareMathOperator{\im}{Im}
\DeclareMathOperator{\fun}{Fun}

\begin{document}

\title{
A simplicial groupoid for plethysm
}
\author{Alex Cebrian}
\date{}
\maketitle

\begin{abstract}
We give a simple combinatorial model for plethysm. Precisely, the bialgebra dual to plethystic substitution is realised as the homotopy cardinality of the incidence bialgebra of an explicit simplicial groupoid, obtained from surjections by a construction reminiscent of the Waldhausen S and the Quillen Q-construction.
\end{abstract}

{\let\thefootnote\relax\footnote{This work has received support from grant number MTM2016-80439-P of Spain.}}

\section*{Introduction}

Plethysm is a substitution operation in the ring of power series in infinitely many variables. It was introduced
by P\'olya  \cite{polya1937} in unlabelled enumeration theory in combinatorics motivated as a series analogue of the wreath product of permutation groups. Another notion of  plethysm was defined by D. E. Littlewood  \cite{Littlewood305}  in the context of symmetric functions and representation theory of the general linear groups \cite{Macdonald}. It appears also in algebraic topology, in connection with $\l$-rings \cite{Borger-Wieland} and power operations in cohomology \cite{bauer}.  The two notions of plethysm are closely related, as described in \cite{StanleyII} and \cite{Bergeron:2}. The present paper is concerned with P\'olya's notion.

The combinatorics involved in plethysm  is rather intricate. In the process of categorification of generating functions, Joyal \cite{Joyal:1981} presented a combinatorial model for the plethystic substitution of cycle index series arising from species. Specifically, he proved that composition of species corresponds to plethystic substitution of their cycle index series.
However, a fully combinatorial construction was only given a few years later by Nava and Rota \cite{Nava-Rota}.
They developed the notion of partitional, a functor from the groupoid of partitions to the category of finite sets, and showed that a suitable notion of composition of partitionals yields plethystic substitution of their generating functions, in analogy with  composition of species and composition of their exponential generating functions. 
A variation of this combinatorial interpretation was given shortly after by Bergeron  \cite{Bergeron}, who instead of partitionals considered  permutationals, functors from the groupoid of permutations to the category of finite sets. This approach  is nicely related to the theory of species and their cycle index series through an adjunction.

 The present work provides a very different construction modeled on Joyal's construction of the Fa\`a di Bruno bialgebra from the category of surjections \cite{Joyal:1981}, which in a modern reformulation \cite{GKT:DSIAMI-2} states that ``the homotopy cardinality of the incidence bialgebra of the fat nerve of the category of surjections, $N\sur\colon \simplexcategory^{\text{op}}\rightarrow \infgrpd$, is isomorphic to the Fa\`a di Bruno bialgebra'' (Section \ref{proof}). Thus, plethystic substitution is exhibited as a convolution tensor product obtained from an explicit simplicial groupoid, $T\sur\colon \simplexcategory^{\text{op}}\rightarrow \infgrpd$, by standard general constructions: incidence coalgebras and homotopy cardinality.   This simplicial groupoid arises from the category of finite sets and surjections $\sur$ as its $T$-\emph{construction}, a categorical construction that we introduce in {Section \ref{T}}. 
 
Let us briefly recall the definition of plethystic substitution, explain the
appearance of $T\sur$ and state our main theorem.

 \vspace{8pt}
 \noindent\textbf{Plethysm.} Given a formal power series 
$$F(x_1,x_2,\dots)=\sum_{(\mu_1,\mu_2,\dots)}f_{(\mu_1,\mu_2,\dots)}\frac{x_1^{\mu_1}x_2^{\mu_2}\dots}{1!^{\mu_1}\mu_1!\;2!^{\mu_2}\mu_2!\cdots} \in \Q\pow{x_1,x_2,\dots},$$
where the sum ranges over over all nonzero sequences $(\mu_1,\mu_2,\dots)$ of nonnegative integers with a finite number of nonzero entries,
 and given a second such formal power series $G(x_1,x_2,\dots)$, the \emph{plethystic substitution} of $F$ into $G$, denoted $G\oast F$, is defined as
 $$(G\oast F)(x_1,x_2,\dots):=G(F_1,F_2,\dots),$$
 where
\begin{equation*}\label{fk}
F_k(x_1,x_2,\dots):=F(x_k,x_{2k},\dots).
\end{equation*}
Observe that $(G\oast F)(x_1,0,\dots)=G(F(x_1,0,\dots),0,\dots)$, the usual composition of one-variable power series. 
In Section \ref{pleth} the plethystic version of the Fa\`a di Bruno bialgebra, called \emph{plethystic bialgebra} and denoted  $\mathcal{P}$, will be studied .  As an algebra $\mathcal{P}$ is free on the linear functionals $A_{\l}(F)=F_{\l}$, and its comultiplication is dual to plethystic substitution. In particular we will derive an explicit formula (Proposition \ref{prop:pletform}) for this comultiplication which  can be viewed as a generalization of the Bell polynomial expression for extracting the coefficients of ordinary composition of one-variable power series.


\vspace{8pt}
 \noindent\textbf{The simplicial groupoid $T\sur$.} Objects in $T_1\sur$ and $T_2\sur$ (1 and 2-simplices of $T\sur$) are, respectively, diagrams
 
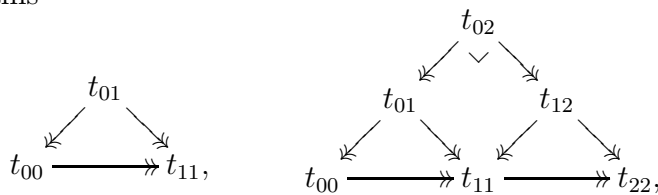
\vspace{-20pt}
 \[
 \xymatrix@R=2.8ex@C=2ex{& &\\
 &{t_{01}}\ar@{->>}[ld] \ar@{->>}[rd]&\\
 	     {t_{00}} \ar@{->>}[rr]& &{t_{11}},} \hspace{30pt}    
\xymatrix@R=2.8ex@C=2ex{& & {t_{02}}\dpullbackS \ar@{->>}[ld] \ar@{->>}[rd]& &\\                             
			&{t_{01}}\ar@{->>}[ld] \ar@{->>}[rd]& &{t_{12}}\ar@{->>}[ld] \ar@{->>}[rd] & \\
	     {t_{00}} \ar@{->>}[rr]& &{t_{11}} \ar@{->>}[rr]& &{t_{22}},}        
\]
where the $t_{ij}$ are finite sets and the arrows surjections. Morphisms of such shapes are levelwise bijections $t_{ij}\xrightarrow{\sim} t_{ij}'$ compatible with the diagram. In general $T_n\sur$ is an analogous pyramid, with $t_{0n}$ in the peak, all of whose squares are pullbacks of sets. The  face maps $d_i$ remove all the sets containing an $i$ index, and the  degeneracy maps $s_i$ repeat the $i$th diagonals, as will be detailed in Section \ref{T}. 

 Diagrams whose last set is
  singleton are called \emph{connected}.
  The main point of $T\sur$ is that
  the connected objects in $T_1\sur$
  parametrise precisely the summation
  of the series, and, rather strikingly,
  the connected objects in $T_2\sur$ encode
  all the combinatorics of plethystic
  substitution, as we shall see. Including also the non-connected objects
  is essential for having a simplicial
  object.
 
Note that the nerve of the category of surjections, $N\sur$, is contained in $T\sur$ as the simplices whose left-down arrows are identities. This reflects the already mentioned fact that plethysm restricts to ordinary composition in the first variable.

\vspace{8pt}
\noindent\textbf{Proposition (cf. \ref{TmonoSegal}).} \textit{The simplicial groupoid $T\sur$ is a CULF monoidal Segal space.}
\vspace{4pt}

The meaning of this words will be explained in Section \ref{preliminaries}. The upshot is that,
in particular, $T\sur$ is a \emph{monoidal decomposition space} in the sense of \cite{GKT:DSIAMI-1}, and therefore, by general principles,  $\infgrpd_{/T_1\sur}$, the slice category of groupoids over $T_1\sur$, becomes a bialgebra at the objective level, meaning that comultiplication and multiplication are linear functors. This bialgebra is called the incidence bialgebra, and its comultiplication is the functor $(d_2,d_0)_!\circ d_1^{\ast}$ (where upperstar is pullback and lowershriek is postcomposition) given by the span of face maps
$$T_1\sur\xleftarrow{\;\;d_1\;  \;} T_2\sur \xrightarrow{(d_2,d_0)\;}T_1\sur\times T_1\sur.$$

Furthermore, $T\sur$ is locally finite. This is the condition ensuring we can take homotopy cardinality of its incidence bialgebra $\infgrpd_{/T_1\sur}$ to get a bialgebra structure on $\Q_{\pi_0T_1\sur}$ (cf. \cite{GKT:DSIAMI-2}).

\vspace{8pt}
 \noindent\textbf{Main Theorem (cf.  \ref{main}).} \textit{The homotopy cardinality of the incidence bialgebra of $T\sur$ is isomorphic to $\mathcal{P}$.}

\vspace{4pt}

In this way plethystic substitution is derived from general principles,
and the relationship between the combinatorial construction
and the algebra is remarkably clean: that of taking homotopy cardinality, whereby all symmetry factors
come out right automatically. 

The notions of Segal space, incidence bialgebra and homotopy cardinality will be reviewed in Section \ref{HoCard} following G\'alvez--Kock--Tonks \cite{GKT:DSIAMI-1, GKT:DSIAMI-2}. The $T$-construction is explained in Section \ref{T}, and the plethystic bialgebra is studied in Section \ref{pleth}. Section \ref{proof} is devoted to the proof of the main theorem.  Section \ref{partitions} makes a synthetic translation of some partition-theoretic constructions used in the Nava--Rota combinatorial interpretation to the language of surjections and $T\sur$. Finally, in Section \ref{FdB} we derive a Fa\`a di Bruno  formula for the ``Green function'' of all surjections, in the style of \cite{GKT:FdB}. 
 
\subsection*{Acknowledgements} The author would like to thank Joachim Kock for his advice, help and support throughout this project, and Andr\'e Joyal and François Bergeron for their feedback.

 
\section{Segal spaces, incidence bialgebras and homotopy cardinality}\label{preliminaries}\label{HoCard}

 We denote by $\infgrpd$ the category of groupoids.  Throughout  the article, pullbacks and fibres of groupoids refer to homotopy pullbacks and homotopy fibres. A brief introduction to the homotopy approach to groupoids in combinatorics can be found in \cite[\S 3]{GKT:FdB}. A simplicial groupoid  $X\colon \simplexcategory^{\text{op}}\longrightarrow \infgrpd$  is a \emph{Segal space} \cite[\S 2.9, Lemma 2.10]{GKT:DSIAMI-1} if the following square is a pullback for all $n>0$:
\begin{equation}\label{segal}
\vcenter{\xymatrix{X_{n+1} \ar[r]^{d_0} \ar[d]_{d_{n+1}} \drpullback& X_n \ar[d]^{d_n}\\
	      X_n \ar[r]_{d_0} & X_{n-1}.}}
\end{equation}

Segal spaces arise prominently through the fat nerve construction: the fat nerve of a category $\mathcal{C}$ is the simplicial groupoid $X=N\mathcal{C}$ with $X_n=\fun ([n],\C)^{\simeq}$, the groupoid of functors $[n] \to \C$. In this case the pullbacks above are strict, so that all the simplices are strictly determined by $X_0$ and $X_1$, respectively the objects and arrows of $\C$, and the inner face maps are given by composition of arrows in $\C$. In the general case $X_n$, is determined from $X_0$ and $X_1$ only up to equivalence, but one may still think of it as a ``category'' object whose composition is defined only up to equivalence.

Let $X$ be a simplicial groupoid.
The spans
$$X_1\xleftarrow{\;\;d_1\;  \;} X_2 \xrightarrow{(d_2,d_0)\;}X_1\times X_1,    \;\;\;\;\;\;\;\;\;\;\;\;\;\; X_1\xleftarrow{\;\;s_0\; \;} X_0 \xrightarrow{\;\;t\;\;}1, $$

\noindent define two functors 

\begin{center}
\begin{tabular}{rllccrll}
$\Delta \colon \infgrpd_{/X_1}$& $\longrightarrow$ &$\mathbf{Grpd}_{/X_1\times X_1}$ && & $\epsilon \colon \infgrpd_{/X_1}$& $\longrightarrow$ &$\infgrpd$\\
$S\xrightarrow{s} X_1$& $\longmapsto$ & $(d_2,d_0)_!\circ d_1^{\ast}(s)$, && &$S\xrightarrow{s} X_1$& $\longmapsto$ & $t_!\circ s_0^{\ast}(s)$ .
\end{tabular}
\end{center} 
Recall that upperstar is pullback and lowershriek is postcomposition. This is the general way in which spans interpret homotopy linear algebra \cite{GKT:HLA}.

Segal spaces are a particular case of decomposition spaces \cite[Proposition 3.7]{GKT:DSIAMI-1}, simplicial groupoids with the property that the functor $\Delta$ is coassociative with the functor $\epsilon$ as counit (up to homotopy). In this case $\Delta$ and $\epsilon$ endow $\infgrpd_{/X_1}$ with a coalgebra structure \cite[\S 5]{GKT:DSIAMI-1}   called the \emph{incidence coalgebra} of $X$. Note that in the special case where $X$ is the nerve of a poset, this
construction becomes the classical incidence coalgebra construction after
taking cardinality, as we shall do shortly.

The morphisms of decomposition spaces that induce coalgebra homomorphisms are the so-called \emph{CULF} functors \cite[\S 4]{GKT:DSIAMI-1}, standing for conservative and unique-lifting-of-factorisations. A Segal space $X$ is \emph{CULF monoidal} if it is a monoid object in the monoidal category $(\mathbf{Dcmp}^{\text{CULF}},\times,1)$ of decomposition spaces and CULF functors \cite[\S 9]{GKT:DSIAMI-1}. More concretely, it is CULF monoidal
 if there is a product $X_n\times X_n\rightarrow X_n$ for each $n$, compatible with the degeneracy and face maps, and such that for all $n$ the squares
\begin{equation}
\label{monoidalpullback}
\vcenter{\xymatrix@R=6ex@C=8ex{X_n\times X_n \ar[r]^{g\times g} \ar[d] \drpullback& X_1\times X_1 \ar[d]\\
	      X_n \ar[r]_{g} & X_1,}}
\end{equation}
where $g$ is induced by the unique endpoint-preserving map $[1]\rightarrow [n]$, are pullbacks \cite[\S  4]{GKT:DSIAMI-1}. For example the fat nerve of a monoidal extensive category is a CULF monoidal Segal space. Recall that a category $\C$ is monoidal extensive if it is monoidal $(\C,+,0)$ and the natural functors $\C_{/A}\times \C_{/B}\rightarrow \C_{/A+B}$  and $\C_{/0}\rightarrow 1$ are equivalences.
 
If $X$ is CULF monoidal then the resulting coalgebra is in fact a bialgebra \cite[\S 9]{GKT:DSIAMI-1}, with product given by
\begin{center}
\begin{tabular}{rllccrll}
$\grpdSprod\colon \infgrpd_{/X_1}\otimes \infgrpd_{/X_1}$&$\xrightarrow{\; \sim \;}$& $\infgrpd_{/X_1\times X_1}$  &$\xrightarrow{\;\; +_! \;\;}$ &$\infgrpd_{/X_1}$ \\
$(G\rightarrow X_1)\otimes (H\rightarrow X_1)$& $\longmapsto$ & $G\times H\rightarrow X_1\times X_1$&$ \longmapsto $& $ G\times H\rightarrow X_1$.
\end{tabular}
\end{center} 
Briefly, a product in $X_n$ compatible with the simplicial structure endows $X$ with a product, but in order to be compatible with the coproduct it has to satisfy the diagram \eqref{monoidalpullback} (i.e. it has to be a CULF functor).

A groupoid $X$ is \emph{finite} if $\pi_0(X)$ is a finite set and $\pi_1(x)=\aut(x)$ is a finite group for every point $x$. If only the latter is satisfied then it is called \emph{locally finite}. A morphism of groupoids is called finite when all its fibres are finite. The \emph{homotopy cardinality} \cite[\S 3]{GKT:HLA} of a finite goupoid $X$ is defined as 
$$|X|:=\sum_{x\in \pi_0X}\frac{1}{|\aut( x)|}\in \Q,$$
and the homotopy cardinality of a finite map of groupoids $A\xrightarrow{p}B$ is 
$$|p|:=\sum_{b\in \pi_0B}\frac{|A_b|}{|\aut(b)|}\delta_b,$$
in the completion of $\Q_{\pi_0B}$, the vector space spanned by $\pi_0B$. In this sum $A_b$ is the homotopy fibre and $\delta_b$ is a formal symbol representing the isomorphism class of $b$. A simple computation shows that $|1\xrightarrow{\name{b}}B|=\delta_b$.

 A Segal space $X$ is \emph{locally finite} \cite[\S 7]{GKT:DSIAMI-2} if $X_1$ is a locally finite groupoid and both $s_0\colon X_0\rightarrow X_1$ and $d_1\colon X_2\rightarrow X_1$ are finite maps. In this case one can take homotopy cardinality to get a comultiplication
\begin{center}
\begin{tabular}{rll}
$\Delta \colon \Q_{\pi_0X_1}$& $\longrightarrow$ &$ \Q_{\pi_0X_1}\otimes  \Q_{\pi_0X_1}$\\
$|S\xrightarrow{s} X_1|$& $\longmapsto$ & $|(d_2,d_0)_!\circ d_1^{\ast}(s)|$
\end{tabular}
\end{center} 
and similarly for $\epsilon$ (cf. \cite[\S 7]{GKT:DSIAMI-2}). Moreover, if $X$ is CULF monoidal then $\Q_{\pi_0X_1}$ acquires a bialgebra structure with the product $\cdot=|\grpdSprod|$. In particular, if we denote by $+$ the monoidal product in $X$, then
$\delta_{a}\cdot \delta_{b}=\delta_{a+b}$ for any $|1\xrightarrow{\name{a}}X_1|$ and $|1\xrightarrow{\name{b}}X_1|$.



\section{The $T$-construction}
\label{T}
The $T$-construction for surjections was given already in the
introduction. For the sake of proving Proposition \ref{TmonoSegal} below, we now give it a more formal treatment, inspired by Lurie's account \cite[\S 1.2.2]{Lurie:HA} of the Waldhausen S-construction. The use of transversal complexes (which we introduce here) instead of ``gap complexes'', on the other hand, is reminiscent of Quillen's Q-construction, which is the ``twisted arrow category'' (edgewise subdivision) of S.

The $T$-construction can be applied to any category which possesses a class of \emph{distinguished squares} satisfying  the following axioms inspired by the properties of pullbacks:
\begin{enumerate}[i)]
\item The identity squares are distinguished and the class of distinguished squares is closed under isomorphisms. Equivalently
\[
\vcenter{\xymatrix{ \ar[d] \ar[r]^{\sim} &   \ar[d] \\
               \ar[r]_{\sim}& } } \text{ is distinguished}.  
\]

\item Given
\[
\vcenter{\xymatrix{ \ar[d] \ar[r] &   \ar[d] \ar[r] &   \ar[d]  \\
               \ar[r]& \ar[r]&} },   
\]
if the right and the left squares are distinguished then the outer rectangle is distinguished.
\item For any two maps $b\xrightarrow{p} a, c\xrightarrow{q} a$  there are maps $d\xrightarrow{p'} c,d\xrightarrow{q'} b$ making a distinguished square. Moreover, for any other two maps $ d'\xrightarrow{p''} c, d'\xrightarrow{q''} b$ making a distinguished square there is a unique isomorphism $d'\xrightarrow{\phi}d$ making the diagram commute,
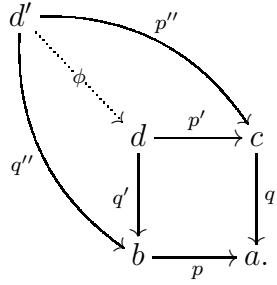
\[
\xymatrix@C=6ex@R=6ex{d' \ar@/_1.3pc/[ddr]_{q''} \ar@/^1.3pc/[drr]^{p''} \ar@{.>}[dr]|-{\phi}\\
&d\ar[d]_{q'} \ar[r]^{p'}& c\ar[d]^q \\
&b \ar[r]_{p} &a.}
\]
 
\end{enumerate}

For instance in any category with pullbacks these form a class of distinguished squares. Moreover, in any subcategory of a category with pullbacks whose arrows are stable under pullbacks these form again a class of distinguished squares. This is our motivating example, since surjections are stable under pullbacks in the category of sets, although the category of finite sets and surjections does not
have pullbacks. Distinguished squares will be indicated with the same symbol as pullbacks. 
 
Let $I$ be a linearly ordered set. Consider the category $\tw(I)$, whose objects are pairs $i\le j$ in $I$ and whose morphisms are relations $(i,j)\le(i',j')$ whenever $i'\le i$ and $j\le j'$ or whenever $i=j\le i'=j'$. This construction can be viewed as the twisted arrow category of $I$ together with arrows between the identities of $I$. 

For example, for $I=[n]$ the objects and arrows of $\tw([n])$ can be pictured as (picturing $n=3$)
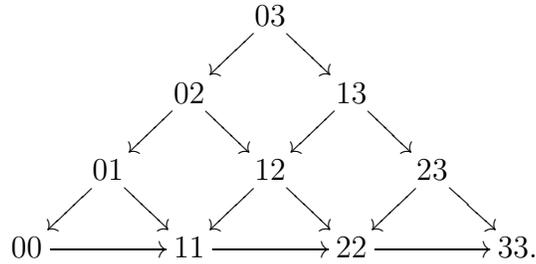
\[
\xymatrix@C=2.5ex@R=3ex{& & &{03} \ar[ld] \ar[rd]& &&\\
		& & {02} \ar[ld] \ar[rd]& & {13} \ar[ld] \ar[rd]& &\\                             
			&{01}\ar[ld] \ar[rd]& &{12}\ar[ld] \ar[rd] & & 23 \ar[ld] \ar[rd]&\\
	     {00} \ar[rr]& &{11} \ar[rr]& &{22} \ar[rr]& &{33}.&}        
\]

The set of categories $\tw([n])$ for all $n$ form a cosimplicial object $\Delta \rightarrow \text{Cat}$ given by $[n]\mapsto \tw([n])$. The face map $d_k\colon  \tw([n{-}1])\rightarrow \tw([n])$ is the obvious induced map $d_k(i,j)=(d_k(i),d_k(j))$, and similarly for the degeneracy maps.

Let $\mathcal{C}$ be a category  with distinguished squares. A functor $F\colon \tw(I)\rightarrow \mathcal{C}$ is called \emph{transversal complex}  if for every $i\le j\le k \le l$ the associated diagram
\begin{equation}\label{Tpullback}
\vcenter{\xymatrix{F(i,l) \ar[d] \ar[r] \drpullback& F(j,l)  \ar[d] \\
              F(i,k) \ar[r]&F(j,k) }  }     
\end{equation}
is a distinguished square, as indicated in the picture. The word transversal complex will be justified in Section \ref{partitions}. Let $\trans(I,\mathcal{C})$ be the full subgroupoid of $\fun (\tw(I), \mathcal{C})^{\simeq}$ containing only the transversal complexes. Then the assignment
$$[n]\longmapsto \trans ([n],\mathcal{C})$$
defines a simplicial groupoid $T\mathcal{C}\colon \simplexcategory^{\text{op}} \rightarrow \infgrpd$.  The groupoid $T_n\C=\trans([n],\C)$ has as objects  diagrams in $\C$ (picturing $n=4$)

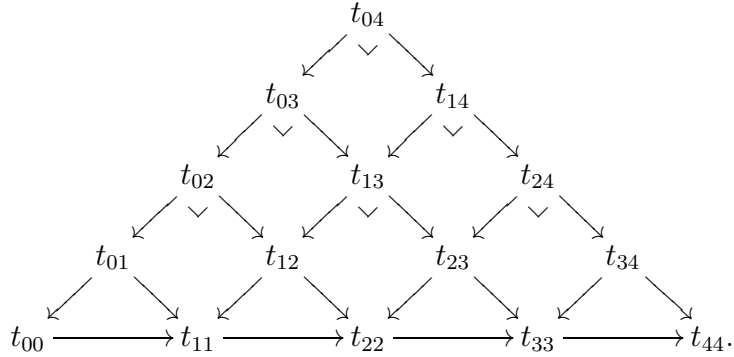
\[
\xymatrix@C=2.5ex@R=3ex{    &&&& t_{04} \ar[ld] \ar[rd] \dpullbackS&&&&\\
		& & & t_{03} \ar[ld] \ar[rd] \dpullbackS& & t_{14}\ar[ld] \ar[rd]\dpullbackS&&&\\
		& & t_{02} \ar[ld] \ar[rd]\dpullbackS& & t_{13} \ar[ld] \ar[rd]\dpullbackS& & t_{24}\ar[ld] \ar[rd]\dpullbackS&&\\                             
			&t_{01}\ar[ld] \ar[rd]& & t_{12}\ar[ld] \ar[rd] & & t_{23} \ar[ld] \ar[rd]& & t_{34}\ar[ld] \ar[rd]&\\
	      t_{00} \ar[rr]& &t_{11} \ar[rr]& &t_{22} \ar[rr]& &t_{33} \ar[rr]& &t_{44}.}        
\]
 The morphisms of such diagrams are levelwise isomorphisms $t_{ij}\xrightarrow{\sim}t'_{ij}$ making the diagram commute. In particular $T_0\C=\C^{\simeq}$. The face map $d_i$ removes all the objects containing an $i$ index. The degeneracy map $s_i$ repeats the $i$th diagonals.
For example
\[
s_1\left(\vcenter{\xymatrix@C=2.5ex@R=3.5ex{& & t_{02} \ar[ld] \ar[rd]\dpullback& & \\                             
			&t_{01}\ar[ld] \ar[rd]& & t_{12}\ar[ld] \ar[rd]&\\
	      t_{00} \ar[rr]& &t_{11} \ar[rr]& &t_{22},}}\right) =
	      \vcenter{\xymatrix@C=2.5ex@R=3.5ex{& & & t_{02} \ar[ld] \ar[rd] \dpullback&&&\\
		& & t_{01} \ar@{=}[ld] \ar[rd]\dpullback& & t_{12} \ar[ld] \ar@{=}[rd]\dpullback&&\\                             
			&t_{01}\ar[ld] \ar[rd]& & t_{11}\ar@{=}[ld] \ar@{=}[rd] & & t_{12} \ar[ld] \ar[rd]&\\
	      t_{00} \ar[rr]& &t_{11} \ar@{=}[rr]& &t_{11} \ar[rr]& &t_{22}.}}
\]
\begin{remark} The Quillen Q-construction of an abelian category $\mathcal{A}$, denoted $Q\mathcal{A}$, can be described in a similar way. It is the simplicial groupoid such that $Q_n\mathcal{A}$ is the full subgroupoid of  $\fun (\text{Tw}([n]), \mathcal{A})^{\simeq}$ consisting of functors $F$ satisfying the same pullback condition as the transversal complexes \eqref{Tpullback} and the additional conditions that for all $i\le j \le k$ the map $F(i,k)\twoheadrightarrow F(i,j)$ is an epimorphism and the map $F(i,k)\rightarrowtail F(j,k)$ is a monomorphism. Thus, an object of  $Q_2\mathcal{A}$ is essentially a diagram
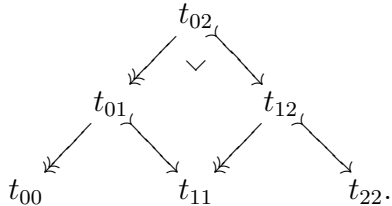
\[
\vcenter{\xymatrix@C=2.5ex@R=3.5ex{& & t_{02} \ar@{->>}[ld] \ar@{>->}[rd]\dpullback& & \\                             
			&t_{01}\ar@{->>}[ld] \ar@{>->}[rd]& & t_{12}\ar@{->>}[ld] \ar@{>->}[rd]&\\
	      t_{00} & &t_{11} & &t_{22}.}}
\]
The main difference between $T$ and $Q$ are the horizontal arrows $a_{ii}\rightarrow a_{jj}$ which appear in $T$ but not in $Q$, coming from the additional arrows $(i,i)\rightarrow (j,j)$ of $\tw([n])$ compared to $\text{Tw}([n])$. For the significance of these extra arrows see Section \ref{partitions}.
\end{remark}

We proceed to show that $T\C$ is a Segal space. 
 We will make use of the following standard result. 
\begin{lemma} Let $F\colon \mathcal{A}\rightarrow \mathcal{B}$ be a functor between two categories and let $\mathcal{C}$ be another category. Then if $F$ is injective on objects the induced functor $\fun (\mathcal{B},\mathcal{C})\rightarrow \fun(\mathcal{A},\mathcal{C})$ is an isofibration.

\end{lemma}
\begin{proposition}\label{T-Segal} Let $\C$ be a category with distinguished squares. Then the simplicial groupoid $T\C$ is a Segal space.
\begin{proof} By the previous lemma, since the face maps $d_i\colon \tw([n{-}1])\rightarrow  \tw([n])$ are injective on objects, the face maps $$d_i\colon \fun(\tw([n]),\C)\rightarrow \fun(\tw([n{-}1]),\C)$$ are isofibrations. But because $T_n\C$ is a full subgroupoid of $\fun(\tw([n]),\C)^{\simeq}$ and is closed under isomorphisms it follows that the face maps of $T\C$ are isofibrations.
As a consequence  it is sufficient to see (see diagram \eqref{segal}) that $T_{n+1}\C$ is equivalent to the strict pullback
\begin{equation}
\vcenter{\xymatrix{ \ar[r] \ar[d] \drpullback&T_n\C \ar[d]^{d_n}\\
	      T_n\C \ar[r]_{d_0} & T_{n-1}\C.}}
\end{equation}
But this is indeed the case: the objects of this pullback are pairs in $T_n\C\times T_n\C$ coinciding at the last and first face respectively. That is, pairs
\[
\left(\vcenter{\xymatrix@C=3.5ex@R=3.5ex{& & \ast \ar[ld] \ar[rd]\dpullback& & \\                             
			&\ast\ar[ld] \ar[rd]& & \ast \ar@{-->}[ld] \ar@{-->}[rd]&\\
	      \ast \ar[rr]& &\ast \ar@{-->}[rr]& &\ast,}}
	    \vcenter{\xymatrix@C=3.5ex@R=3.5ex{& & \ast \ar[ld] \ar[rd]\dpullback& & \\                             
			&\ast\ar@{-->}[ld] \ar@{-->}[rd]& & \ast\ar[ld] \ar[rd]&\\
	     \ast \ar@{-->}[rr]& &\ast \ar[rr]& &\ast}}\right)
\]
such that the dashed regions are equal. This gives a groupoid whose objects are 
 \[
	      \vcenter{\xymatrix@C=3.5ex@R=3.5ex{
		& & \ast \ar[ld] \ar[rd]\dpullback& & \ast \ar[ld] \ar[rd]\dpullback&&\\                             
			&\ast\ar[ld] \ar[rd]& &\ast\ar[ld] \ar[rd] & & \ast \ar[ld] \ar[rd]&\\
	      \ast \ar[rr]& &\ast \ar[rr]& &\ast \ar[rr]& &\ast}}
\]
and morphisms are levelwise isomorphisms making the diagram commute. But this is equivalent to $T_{n+1}\C$, since the missing apex can be uniquely filled with a distinguished square.
\end{proof}
\end{proposition}
The following result characterizes the categories whose $T$-construction is CULF mo\-noidal. 
\begin{proposition}\label{monodec} Let $(\C,+,0)$ be a monoidal category with distinguished squares. Then $T\C$ is a CULF monoidal Segal space if and only if the following conditions hold.
\begin{enumerate}[i)]
\item $(\C,+,0)$ is monoidal extensive, 
\item Given $\Omega,\Gamma, \Lambda$ commutative squares in $\C$ such that $\Omega=\Gamma+\Lambda$, then 
$$\Omega \text{ is distinguished  }\Longleftrightarrow \Gamma \text{ and } \Lambda \text{ are distinguished}.$$
\end{enumerate}
\begin{proof} Condition i) is the necessary and sufficient condition to be able to sum arrows and commutative diagrams in $\C$, as we know from the nerve of a monoidal extensive category. However we need sums of distinguished squares to be distinguished squares in order for $T_n\C\times T_n\C\rightarrow T_n\C$ to be well-defined. This gives the $\Leftarrow$ of ii). Finally, we have to impose $\Rightarrow$ of ii) to ensure that the square \eqref{monoidalpullback} is a pullback.
\end{proof}
\end{proposition}

As explained in the introduction, the central object of this article is $T\sur$, the $T$-construction of the category $\sur$ of finite sets and surjections. Let us now see that $T\sur$ meets all the requirements for the main theorem to be stated.
\begin{lemma} The category $\sur$
\begin{enumerate}[i)]
\item has a class of distinguished squares,
\item is monoidal extensive with disjoint union $(+)$ and empty set as monoidal structure,
\item satisfies ii) of Proposition \ref{monodec},
\end{enumerate}
\begin{proof} We declare the distinguished squares to be the commutative squares in $\sur$ that are pullbacks in the category of sets (note that $\sur$ itself does not have pullbacks. See Section \ref{partitions} for the combinatorial relevance of these subtleties). For ii) observe that taking disjoint union clearly gives an equivalence $\sur_{/A}\times \sur_{/B} \simeq \sur_{/A+B}$. It is the restriction to surjections of the monoidal structure of finite sets and their coproduct. Now, since the pullback in sets is in fact the disjoint union of the product of the fibres we obtain iii). 
\end{proof}
\end{lemma}
In view of this lemma and Propositions \ref{T-Segal} and \ref{monodec} we obtain the following result. 
\begin{proposition}\label{TmonoSegal}$T\sur$ is a CULF monoidal Segal space.
\end{proposition} 
This implies that $T\sur$ has an associated incidence bialgebra (\S \ref{preliminaries}).
Now, observe that $\sur^{\simeq}$ is locally finite, so that $T_1\sur$ is locally finite and $s_0\colon T_0\sur\rightarrow T_1\sur$ is finite. Moreover, every arrow of $\sur$ admits, up to isomorphism, a finite number of 2-step factorizations, therefore  $d_1\colon T_2\sur \rightarrow T_1\sur$ is also finite. This means that $T\sur$ is locally finite in the sense of \cite{GKT:DSIAMI-2}. As a consequence we can take homotopy cardinality of the incidence bialgebra of $T\sur$.


\section{Plethystic bialgebra}
\label{pleth}
It was stated in the introduction how plethystic substitution works. In this section we introduce the bialgebra $\mathcal{P}$, whose comultiplication is dual to plethystic substitution, and derive a formula \eqref{eq:pletform} for extracting the comultiplication of the elements of its basis. The following notation is used.

\begin{multicols}{2}
\begin{itemize}
\item $\x=(x_1,x_2,\dots)$,
\item $\Lambda$: the set of all infinite vectors with a finite number of nonzero entries,
\item $\Lambda\ni\l=(\l_1,\l_2\dots)$,
\item $|\l|=\sum_k \l_k$,
\item $\x^{\l}=x_1^{\l_1}x_2^{\l_2}\cdots$,
\item $\autiv (\l)=1!^{\l_1}\l_1!\cdot 2!^{\l_2}\l_2!\cdots$,
\item $\l+\mu$ is coordinate-wise sum.
\end{itemize}
\end{multicols}

First of all, we define the \emph{$n$-th Verschiebung operator} $V^n$ as

\begin{eqnarray*}
    (V^n\l)_i=\left \{ \begin{array}{ll} \l_{i/n} & \mbox{if} \; n\mid i \\
        0 & \mbox{otherwise,} 
    \end{array} \right. \text{ for }  \l=(\l_1,\l_2\dots).
\end{eqnarray*}

\noindent For example
$V^2(5,9,2,0\dots)=(0,5,0,9,0,2,0\dots).$  It is clear that $V^n$ preserves sums. Note that the $F_k$ defined in the introduction can be expressed as
$$F_k(x_1,x_2,\dots)=F(x_k,x_{2k},\dots)=\sum_{\mu} f_{\mu}\frac{\x^{V^k\mu}}{\autiv (\mu)}.$$

\begin{remark}\label{scalarmultiplication}
Note that $\l$ can be viewed as the isomorphism class of a surjection ${X\twoheadrightarrow B}$ with $\l_k$ fibres of size $k$. With this identification, $\autiv(\l)$ is precisely the cardinal of $\aut(X\twoheadrightarrow B)$ in the groupoid of surjections $\Sg$, whose objects are surjections and whose arrows are pairs of compatible bijections, one for the source and one for the target. Moreover, the Verschiebung operators can also be defined at the objective level of surjections,
$$V^S(X\twoheadrightarrow B):=X\times S\twoheadrightarrow X\twoheadrightarrow B,$$
which is nothing but the scalar multiplication of $X\twoheadrightarrow B$ and $S$ in $\mathbf{Set}_{/B}$ \cite{GKT:HLA}. It is clear that $V^S$ corresponds numerically to $V^{|S|}$.
\end{remark}

For each $\lambda$ define the functional $A_{\l}\in \Q\pow{\x}^{\ast}$ by  $A_{\l}(F)=f_{\l}$. The \emph{plethystic bialgebra} is the free polynomial algebra $\mathcal{P}=\Q[\{A_{\l}\}_{\l}]$ along with the comultiplication dual to plethystic substitution. That is, for each $\lambda$ and $F,G \in \Q\pow{\x}$,

$$\langle \Delta(A_{\l}),F\otimes G\rangle=\langle A_{\l},G\oast F\rangle.$$
The counit is given by $\epsilon (A_{\l})=\langle A_{\l}, x_1\rangle$.

Now, consider a multiset $\pmb{\mu}\in \Lambda^n/\mathfrak{S}_n$ of $n$ infinite vectors. We denote by $\rep(\pmb{\mu})$ the set of automorphisms that map every element to itself. For example if $\pmb{\mu}=\{\alpha,\alpha,\beta,\gamma,\gamma,\gamma\}$ then $\rep(\pmb{\mu})$ has $2!\cdot 1!\cdot 3!$ elements. 
\begin{remark} Observe that $\sum_n \Lambda^n/\mathfrak{S}_n \simeq \pi_0 T_1\sur$. Furthermore, the number of automorphisms of a representative element in $T_1\sur$ of the image of $\pmb{\mu}$ under this bijection is precisely
$$\autiv(\pmb{\mu})= \prod_{\mu\in \pmb{\mu}}\autiv(\mu)\cdot| \rep(\pmb{\mu}) |.$$
\end{remark}

Fix two infinite vectors, $\sigma,\l\in \Lambda$, and a multiset of infinite vectors $\pmb{\mu}\in \Lambda^{|\l|}/\mathfrak{S}_{|\l|}$ . We define the set of $(\l,\pmb{\mu})-$\emph{decompositions} of $\Sg$ as
$$T_{\sigma,\l}^{\pmb{\mu}}:=\left\{p\colon \pmb{\mu}\xrightarrow{\;\;\sim\;\;} \sum_k \{1,\dots,\l_k\}\; |\;\sigma=\sum_{\mu\in\pmb{\mu}}V^{q(\mu)}\mu\right\},$$
where $q$ returns the index of $p(\mu)$ in the sum. A useful way to visualize an element of this set is as a placement of the elements of $\pmb{\mu}$ over a grid with $\l_k$ cells in the $k$th column such that if we apply $V^k$ to the $k$th column and sum the cells the result is $\sigma$. For example, if $\l=(2,0,1,3)$ and $\pmb{\mu}=\{\alpha,\alpha,\beta,\gamma,\gamma,\gamma\}$ the placement

\begin{center}
\begin{tikzpicture}[scale=0.7]
\draw (1,1)--(1,0) -- (0,0) -- (0,1)--(1,1)--(1,2)--(0,2)--(0,1);
\draw (1,0)--(2,0);
\draw (3,1)--(3,0) -- (2,0) -- (2,1)--(3,1);
\draw (3,0)--(4,0)--(4,1)--(3,1)--(3,2)--(4,2)--(4,1)  (3,2)--(3,3)--(4,3)--(4,2);
\node at (0.5,0.5) {$\gamma$}; \node at (0.5,1.5) {$\alpha$}; \node at (0.5,-0.5) {$V^1$};
\node at (2.5,0.5) {$\gamma$}; \node at (1.5,-0.5) {$V^2$}; \node at (2.5,-0.5) {$V^3$};
\node at (3.5,0.5) {$\alpha$}; \node at (3.5,1.5) {$\beta$}; \node at (3.5,2.5) {$\gamma$}; \node at (3.5,-0.5){$V^4$};
\end{tikzpicture}
\end{center}

\noindent belongs to $T_{\sigma,\l}^{\pmb{\mu}}$ if $\sigma=V^1(\gamma+\alpha)+V^3(\gamma)+V^4(\alpha+\beta+\gamma)$. Note that each such placement appears $|\rep (\pmb{\mu})|$ times in $T_{\sigma,\l}^{\pmb{\mu}}$. 

\begin{proposition} \label{prop:pletform} Let $\sigma$ be an infinite vector. Then the comultiplication of $A_{\sigma}$ in $\mathcal{P}$ is given by
\begin{equation}\label{eq:pletform}
\Delta(A_{\sigma})=\sum_{\l}\sum_{\pmb{\mu}} \frac{\autiv(\sigma)\cdot|T_{\sigma,\l}^{\pmb{\mu}}|}{\autiv(\l)\cdot \displaystyle\autiv(\pmb{\mu}) } \prod_{\mu\in \pmb{\mu}} A_{\mu}\otimes A_{\l}.
\end{equation}
\end{proposition}
\begin{remark}Not surprisingly, if $\sigma=(n,0,0,\dots)$ this expression gives the comultiplication of $A_n$ for ordinary composition of one-variable power series.
 Extending this analogy between classical and plethystic we define the polynomials $P_{\sigma,\l}(\{A_{\mu}\}_{\mu})$,
$$\Delta(A_{\sigma})=:\sum_{\l}P_{\sigma,\l}\left(\{A_{\mu}\}_{\mu}\right)\otimes A_{\l},$$
which are the generalization of the Bell polynomials to the plethystic case. Hence in particular $P_{(n,0,\dots),(k,0,\dots)}\left(\{A_{\mu}\}_{\mu}\right)=B_{n,k}\left(\{A_{(i,0,\dots)}\}_i\right)$.
\end{remark}

Before proving Proposition \ref{prop:pletform} we shall need the following two lemmas. For the sake of notation we will work from now on with another basis of $\mathcal{P}$, $\{a_{\l}\}_{\l}$, defined as $a_{\l} :=\frac{A_{\l}}{\autiv({\l})}$.
Let $\z=z_1,z_2,\dots$ be a set of infinitely many formal variables. Consider, in the style of \cite[Remark 2.3]{BFK:0406117}, the map $\Delta \colon \mathcal{P} \pow{\z}\rightarrow (\mathcal{P}\otimes \mathcal{P})\pow{\z}$ given by linearly extending the comultiplication defined above for $\mathcal{P}$. Define the power series
$$A_i(\z)=\sum_{\l}a_{\l}\z^{V^i\l}, \;\;\;\;\;\; i\ge 0.$$
Observe that by definition $\Delta(A_i(\z))=\sum_{\l}\Delta(a_{\l})\z^{V^i\l}$.
The following result is straightforward.
\begin{lemma} \label{Aifi}For any $i,j\ge 0$ and $F,G\in \Q\pow{\x}$
\begin{enumerate}[i)]
\item $\langle A_i(\z),F\rangle =F_i(\z)$,
\item $\langle \Delta(A_{i}(\z)),F\otimes G\rangle = (G\oast F)_i(\z)$,
\item $\langle A_i(\z),F\cdot G\rangle=\langle A_i(\z),F\rangle \cdot \langle A_i(\z),G\rangle$,
\item $\langle A_i(\z)\cdot A_j(\z),F\rangle =\langle A_i(\z),F\rangle \cdot \langle A_j(\z),F\rangle$.
\end{enumerate}
\end{lemma}
\begin{lemma} \label{lemmaCopA}
\begin{equation*}
\Delta(A_1(\z))=\sum_{\l} \left( \prod_i A_i^{\l_i}(\z)\right )\otimes a_{\l}.
\end{equation*}
\begin{proof} Let $F,G\in \Q\pow{\x}$. By ii) of Lemma \ref{Aifi} we have
$$\langle \Delta(A_{1}(\z)),F\otimes G\rangle=(G\oast F) (\z).$$
Now, by definition of plethystic substitution
$$(G\oast F) (\z)=\sum_{\l} \left( \prod_i F_i^{\l_i}(\z)\right )\cdot a_{\l}(G(\z)),$$
but  iv) and i) of Lemma \ref{Aifi} tell us respectively that
$$\left\langle \sum_{\l} \left( \prod_i A_i^{\l_i}(\z)\right ),F\right\rangle=\sum_{\l}  \left( \prod_i \left\langle A_i(\z),F\right\rangle^{\l_i} \right )=\sum_{\l} \left( \prod_i F_i^{\l_i}(\z)\right ).$$
Therefore 
$$(G\oast F) (\z)= \left\langle \sum_{\l} \left( \prod_i A_i^{\l_i}(\z)\right )\otimes a_{\l}, F\otimes G    \right\rangle,$$
as we wanted to see.
\end{proof}
\end{lemma}

\noindent \textit{Proof of Proposition \ref{prop:pletform}.}\\
Define sets
\begin{align*}
T_{\l}:=&\left\{\{\mu_{i,j}\}_{i\ge 1,j\in \{1,\dots, \l_i\}} \right\},\\
T_{\sigma,\l}:=&\left\{\{\mu_{i,j}\}_{i\ge 1,j\in \{1,\dots, \l_i\}} \;|\; \sigma=\sum_i\sum_{j=1}^{\l_i}V^i\mu_{i,j}\right\}.
\end{align*}
We now compute
\begin{align*}
\Delta(A_1(\z))= \sum_{\l} \left( \prod_i A_i^{\l_i}(\z)\right )\otimes a_{\l}=&\sum_{\l}\left( \sum_{\{\mu_{i,j}\}_{i,j} \in T_{\l}}
\prod_i \prod_{j=1}^{\l_i}a_{\mu_{i,j}}\z^{V^i\mu_{i,j}}\right)\otimes a_{\l}=\\ \nonumber
 =& \sum_{\l} \left(\sum_{\sigma} \left(\sum_{\{\mu_{i,j}\}_{i,j} \in T_{\sigma,\l}}
\prod_i \prod_{j=1}^{\l_i}a_{\mu_{i,j}}\right)\z^{\sigma}\right)\otimes a_{\l}=\\ \nonumber
=& \sum_{\sigma} \sum_{\l} \left(\sum_{\{\mu_{i,j}\}_{i,j} \in T_{\sigma,\l}}
\prod_i \prod_{j=1}^{\l_i}a_{\mu_{i,j}}\right)\z^{\sigma}\otimes a_{\l}.\\
\end{align*}
But on the other hand $\Delta(A_1(\z))=\sum_{\sigma}\Delta(a_{\sigma})\z^{\sigma}$. Hence by Lemma  \ref{lemmaCopA} 
\begin{equation}\label{pletform}
\Delta(a_{\sigma})=\sum_{\l}\left(\sum_{\{\mu_{i,j}\}_{i,j}\in T_{\sigma,\l}} \prod_i \prod_{j=1}^{\l_i}a_{\mu_{i,j}}\right)\otimes a_{\l} .
\end{equation}


Notice that sequences $\{\mu_{i,j}\}_{i,j}$ with the same elements ordered in different ways may be $\l$-decompositions of distinct $\sigma$. Therefore the comultiplications of different generators may have terms in common.   
Finally, observe that as a multiset any sequence $\pmb{\mu}=\{\mu_{i,j}\}_{i,j}$ appears precisely $\displaystyle \frac{|T_{\sigma,\l}^{\pmb{\mu}}|}{|\rep(\pmb{\mu})|}$ times in $T_{\sigma,\l}$. This implies that expression \eqref{pletform} is equivalent to 

\begin{equation*}
\Delta(a_{\sigma})=\sum_{\l}\left(\sum_{\pmb{\mu}} \frac{|T_{\sigma,\l}^{\pmb{\mu}}|}{|\rep(\pmb{\mu})|}\prod_{\mu\in\pmb{\mu}}a_{\mu}\right)\otimes a_{\l} .
\end{equation*}
Now, changing again the basis from $\{a_{\l}\}_{\l}$ to $\{A_{\l}\}_{\l}$ we obtain equation \eqref{eq:pletform}.
\qed


\section{Main Theorem}
\label{proof}
We are  now ready to state and prove the main result of this paper. The proof is essentially a question of unpacking the abstract constructions. A pleasant feature is the way in which the subtle symmetry factors come out naturally from the groupoid formalism. 
\vspace{-5pt}
\begin{theorem}\label{main} The homotopy cardinality of the incidence bialgebra of $T\sur$ is isomorphic to $\mathcal{P}$.
\vspace{-3pt}
\begin{proof}
Recall from Section \ref{preliminaries} that the homotopy cardinality of the incidence bialgebra of $T\sur$ is denoted by $\Q_{\pi_0T_1\sur}$. 
We split the proof into three parts. First we define an isomorphism $\Q_{\pi_0T_1\sur}\overset{\theta}{\simeq}\mathcal{P}$ of algebras, next we explore the relation between the Verschiebung operator and $T_2\sur$, and finally we show that $\theta$ preserves the comultiplication. 

\vspace{8pt}
\noindent\textbf{The isomorphism.} We will call \emph{connected} the elements of $T_n\sur$ with a singleton at the $nn$ position. Notice that as a vector space $\Q_{\pi_0T_1\sur}$ is spanned by $\pi_0T_1\sur$, and as a free algebra it is generated by the classes of the connected elements of $T_1\sur$, since every diagram is a sum of connected ones.
The isomorphism class $\delta_{\l}$ of  a connected element 
\[
\vcenter{\xymatrix{& t_{01} \ar@{->>}[ld] \ar@{->>}[rd]  & \\
             t_{00} \ar@{->>}[rr]&  \ar@{}[u]|{\mbox{$\l$}}&1 }}      \in T_1\sur
\]
is given by the infinite vector $\l=(\l_1,\l_2,\dots)$ representing the  class of $t_{01}\twoheadrightarrow t_{00}$. Be aware that the same notation is used for either the connected elements of   $T_1\sur$ and the infinite vectors representing their isomorphism class. This being said, the assignment 
\begin{center}
\begin{tabular}{rll}
$\Q_{\pi_0T_1\sur}$& $\longrightarrow$ &$\mathcal{P}$\\
$\delta_{\l}$& $\longmapsto$ & $A_{\l}$\\
$\delta_{\l+\mu}=\delta_{\l}\delta_{\mu}$&$\longmapsto$ & $A_{\l}A_{\mu}$,
\end{tabular}
\end{center} 
for $\l$ and $\mu$ connected, defines  an isomorphism of algebras. Notice that $\l+\mu$ is the monoidal sum in $T_1\sur$, which does not correspond to the pointwise sum of their corresponding infinite vectors, since it has two connected components.

%

\vspace{8pt}
\noindent\textbf{The Verschiebung operator.} Pick a connected element $t$,
\[
\vcenter{\xymatrix@C=2.5ex@R=3.5ex{& & t_{02} \ar@{->>}[ld] \ar@{->>}[rd]\dpullback& & \\                             
			&t_{01}\ar@{->>}[ld]_{\mu} \ar@{->>}[rd] & & t_{12}\ar@{->>}[ld] \ar@{->>}[rd] &\\
	      t_{00} \ar@{->>}[rr]& &t_{11} \ar@{->>}[rr]& &1}}\in T_2\sur.
\]
For each $r\in t_{11}$, consider the map on the fibres $\mu_r\colon (t_{01})_r\twoheadrightarrow (t_{00})_r$. Since the square is a pullback of sets we have that
 $$(t_{02})_r\longrightarrow (t_{00})_r \simeq (t_{01})_r\times (t_{12})_r \rightarrow (t_{01})_r\xrightarrow{\;\mu_r\;} (t_{00})_r=V^{(t_{12})_r}\mu_r,$$ which  belongs to the isomorphism class of $V^{|(t_{12})_r|}\mu_r$ (see Remark \ref{scalarmultiplication}). But recall that 
 \[
d_1(t)=\vcenter{\xymatrix{& t_{02} \ar@{->>}[ld] \ar@{->>}[rd]  & \\
             t_{00} \ar@{->>}[rr]&  &1 }},
\]
therefore the isomorphism class of $d_1(t)$ is precisely
\begin{equation*}
\sum_{r\in t_{11}} V^{|(t_{12})_r|}\mu_r.
\end{equation*}
In fact, continuing with the scalar multiplication interpretation of Remark \ref{scalarmultiplication}, we could just say that $t_{02}\twoheadrightarrow t_{00}$ is a linear combination in $\bigoplus_r\mathbf{Set}_{/(t_{00})_r}$, namely
$$t_{02}\longrightarrow t_{00}=\sum_r ((t_{02})_r\longrightarrow (t_{00})_r)\cdot (t_{12})_r=\sum_r V^{(t_{12})_r}((t_{02})_r\longrightarrow (t_{00})_r).$$


\noindent\textbf{The comultiplication.} 
We have to show that $\Delta(\delta_{\sigma})$ for $\sigma$ connected yields equation \eqref{eq:pletform}. The first step is a general property of Segal groupoids. Recall that given a locally finite Segal space $X$, the comultiplication $\Delta$ is induced by the span
$$X_1\xleftarrow{\;\;d_1\;  \;} X_2 \xrightarrow{(d_2,d_0)\;}X_1\times X_1$$
by pullback along $d_1$ and postcomposition with $(d_2,d_0)$, and after taking homotopy cardinality
\begin{center}
\begin{tabular}{rll}
$\Delta \colon \Q_{\pi_0X_1}$& $\longrightarrow$ &$ \Q_{\pi_0X_1}\otimes  \Q_{\pi_0X_1}$\\
$|S\xrightarrow{s} X_1|$& $\longmapsto$ & $|(d_2,d_0)_!\circ d_1^{\ast}(s)|.$
\end{tabular}
\end{center} 
Thus, by definition
$$\Delta(\delta_f)=\Delta(|1\xrightarrow{\name{f}}X_1|)=|(X_2)_f\xrightarrow{(d_2,d_0)_!\circ d_1^{\ast}(f)} X_1\times X_1|,$$
so that
$$\Delta(\delta_f)=\sum_{b\in \pi_0X_1}\sum_{a\in \pi_0X_1} \frac{|(X_2)_{f,a,b}|}{|\aut(b)||\aut(a)|}\delta_a\otimes \delta_b,$$
where $(X_2)_{f,a,b}$ stands for homotopy fibre of $f,a$ and $b$ along  $d_1,d_2$ and $d_0$ respectively.
Now, for $X$ a Segal groupoid,  \cite[Lemma 7.10]{GKT:DSIAMI-2} tells us that $(X_{2})_{a,b}$ is discrete and that there is a natural bijection from $(X_{2})_{a,b}$ to $\iso(d_0a,d_1b)$, the set of isomorphisms between $d_0a$ and $d_1b$. Hence, by taking homotopy fibre over $f$ along $d_1$ at both sides we get a natural bijection $(X_{2})_{f,a,b}\simeq \iso(d_0a,d_1b)_f$, and therefore we obtain the following expression for $\Delta(\delta_f)$.
\begin{lemma}
\label{SegalCom} 
Let $X$ be a Segal space. Then for $f$ in $X_1$ we have
$$\Delta(\delta_f)=\sum_{b\in \pi_0X_1}\sum_{a\in \pi_0X_1} \frac{|\iso(d_0a,d_1b)_{f}|}{|\aut(b)||\aut(a)|}\delta_a\otimes \delta_b.$$
\end{lemma}
\noindent In the case $X=T\sur$ we have 
\begin{equation}\label{TSpletform}
\Delta(\delta_{\sigma})=\sum_{\l\in \pi_0 T_1\sur} \sum_{\tau\in \pi_0 T_1\sur}\frac{|\iso(d_0\tau,d_1\l)_{\sigma}|}{|\aut(\l)||\aut(\tau)|}\delta_{\tau}\otimes\delta_{\l}.
\end{equation}
The notation is taken from Proposition \ref{prop:pletform}, meaning that $\delta_{\l}$ corresponds to $A_{\l}$ and $\delta_{\sigma}$ corresponds to $A_{\sigma}$. 

Now, since $\sigma$ is connected also $\l$ must  be connected,
 \[
\vcenter
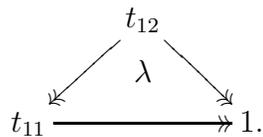
{\xymatrix{ & t_{12} \ar@{->>}[rd]  \ar@{->>}[ld]& \\
             t_{11} \ar@{->>}[rr]&\ar@{}[u]|{\mbox{$\l$}}&1. }}
\]
Now, it is clear that $\iso(d_0\tau,d_1\l)$ is nonempty if and only if  $|d_0\tau|=|d_1\l|$. Without loss of generality we can assume that $d_0(\tau)=d_1(\l)=t_{11}$ and write
 \[
\vcenter
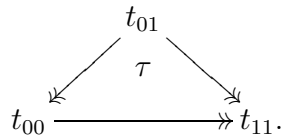
{\xymatrix{ & t_{01} \ar@{->>}[ld]  \ar@{->>}[rd]& \\
             t_{00} \ar@{->>}[rr]&\ar@{}[u]|{\mbox{$\tau$}}&t_{11} .}}
\]
Note that in particular $|t_{11}|=|\l|$. If, as above, we denote by $\mu_r\colon (t_{01})_r\twoheadrightarrow (t_{00})_r$ the fibre surjection of $r\in t_{11}$ and $\pmb{\mu}=\{\mu_r\}_r$, we have that $\delta_{\tau}$ corresponds to $\prod_{\mu \in \pmb{\mu}}A_{\mu}$ and that
$$|\aut(\tau)|= \autiv(\pmb{\mu}).$$
 Hence it only remains to show that $|\iso(d_0\tau,d_1\l)_{\sigma}|=\autiv(\sigma)\cdot |T_{\sigma,\l}^{\pmb{\mu}}|$.
Since  $\iso(d_0\tau,d_1\l)$ is a discrete groupoid, i.e. just a set, it makes sense to consider the subset
$$\{\phi \in \iso(d_0\tau,d_1\l)\; |  \; d_1(\phi)\simeq \sigma\}$$
consisting of those $\phi$ that give an object isomorhpic to $\sigma$ after composition by $\phi$, pullback and $d_1$, as shown in the picture,
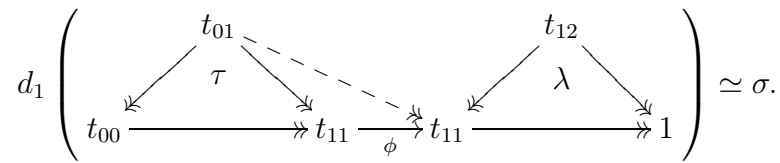
\[
d_1\left(\vcenter{\xymatrix{& t_{01} \ar@{->>}[ld] \ar@{->>}[rd] \ar@{-->>}[rrd]& && t_{12}\ar@{->>}[ld] \ar@{->>}[rd] &\\
             t_{00} \ar@{->>}[rr]&\ar@{}[u]|{\mbox{$\tau$}}&t_{11}\ar[r]_{\phi}&t_{11}\ar@{->>}[rr]&\ar@{}[u]|{\mbox{$\l$}} &1 }}\right)\simeq \sigma.
\]
Note that in this picture the pullback is taken along the dashed arrow. Observe that this condition on $\phi$ can be written as 
$$\sigma=\sum_{r\in t_{11}} V^{|(t_{12})_{\phi(r)}|}\mu_r,$$
where now $\sigma$ and $\mu_r$ represent the corresponding infinite vectors.
This is in bijection with morphisms $\pmb{\mu}\xrightarrow{\sim} \sum \{1,\dots,\l_k\}$ satisfying
$$\sigma=\sum_{\mu\in\pmb{\mu}}V^{q(\mu)}\mu,$$
since summing over $r\in t_{11}$ is equivalent to summing over $\mu \in \pmb{\mu}$ and the same Verschiebung operators appear in both sums. Hence there is a bijection 
$$\{\phi \in \iso(d_0\tau,d_1\l)\; |  \; d_1(\phi)\simeq \sigma\}\simeq T_{\sigma,\l}^{\pmb{\mu}}=\left\{p\colon \pmb{\mu}\xrightarrow{\;\;\sim\;\;} \sum_k \{1,\dots,\l_k\}\; |\;\sigma=\sum_{\mu\in\pmb{\mu}}V^{q(\mu)}\mu\right\}.$$
But the homotopy fibre of $\iso(d_0\tau,d_1\l)$ over $\sigma$ is precisely this subset described above times the set of automorphisms of $\sigma$ in $T_1\sur$, that is
$$\iso(d_0\tau,d_1\l)_{\sigma}\simeq \iso (\sigma,\sigma)\times \{\phi \in \iso(d_0\tau,d_1\l)\; |  \; d_1(\phi)\simeq \sigma\}.$$
Therefore   $|\iso(d_0\tau,d_1\l)_{\sigma}|=\autiv(\sigma)\cdot |T_{\sigma,\l}^{\pmb{\mu}}|$ and equation \eqref{TSpletform} corresponds to equation \eqref{eq:pletform}, as we wanted to see.
\end{proof}
\end{theorem}
\begin{remark} We end this section by deriving the analogous result for one-variable power series. As mentioned in the introduction, the statement reads: the Fa\`a di Bruno bialgebra $\mathcal{F}$ is equivalent to $\Q_{\pi_0\Sg}$, the homotopy cardinality of the incidence bialgebra of the fat nerve $N\sur\colon \simplexcategory^{\text{op}}\rightarrow \infgrpd$ of the category of surjections. First of all, observe that $\mathcal{F}$ is generated by the functionals $A_n$ returning the $n$th coefficient of a power series. The connected elements of $\Sg$ are the surjections with singleton target. Hence,  ${\delta_n=|1\xrightarrow{\name{n\twoheadrightarrow 1}} \Sg|}$ corresponds to $A_n$. Using Lemma \ref{SegalCom} we get
$$\Delta (\delta_n)=\sum_{b:k\twoheadrightarrow 1}\sum_{a:n\twoheadrightarrow k} \frac{|\iso(k,k)_{n\twoheadrightarrow 1}|}{|\aut(b)||\aut(a)|}\delta_a\otimes \delta_k.$$
It is clear that $|\aut(b)|=k!$. Now, in this case any element of $\iso(k,k)$ gives $n\twoheadrightarrow 1$, so that $|\iso(k,k)_{n\twoheadrightarrow 1}|=n!\cdot k!$.  Moreover, $\delta_a=\delta_{n_1}\dots \delta_{n_k}$, where $n_i$ are the fibres of $a\colon n\twoheadrightarrow k$.
Altogether we obtain
$$\Delta (\delta_n)=\sum_{k\twoheadrightarrow 1}\sum_{n\twoheadrightarrow k} \frac{n!}{|\aut(n\twoheadrightarrow k)|}\prod_{i=1}^k \delta_{n_i}\otimes \delta_k,$$
which is easily checked to correspond to the comultiplication of $A_n$,
$$\Delta (A_n)=\sum_{k}\sum_{n_1+\dots +n_k=n} {n\choose n_1,\dots,n_k} \prod_{i=1}^k A_{n_i}\otimes A_k.$$
\end{remark}


\section{Partitions, transversals and $T\sur$}\label{partitions}

As stated in the introduction, the combinatorial interpretation given by Nava--Rota \cite{Nava-Rota} makes use of the groupoid of partitions, whose objects are sets with a partition and whose arrows are block preserving bijections between the underlying sets. The role that partitions of sets play in composition of species is played by a subtle particular case of partition of a partition called transversal in composition of partitionals.
 
The groupoid of partitions is equivalent to the groupoid of surjections $\Sg$. However, the groupoid of surjections has better functorial properties than the groupoid of partitions, mainly because partitions of partitions are pairs of composable surjections. Although the results in this article have been derived independently of the theory of partitionals, it is worth making a brief translation of some crucial partition constructions to the language of surjections.

Consider two partitions $\pi$ and $\tau$ on a set $E$. Let $\pi\colon E\twoheadrightarrow S$ and $\tau\colon E\twoheadrightarrow X$ be their corresponding surjections. Construct the diagram of sets

\begin{equation}\label{partsurj}
\vcenter{\xymatrix@C=6ex@R=6ex{E \ar@{->>}@/_1.3pc/[ddr]_{\pi} \ar@{->>}@/^1.3pc/[drr]^{\tau} \ar@{.>}[dr]|-{\phi}\\
&S \times_I X \drpullback \ar@{->>}[d] \ar@{->>}[r]& X\ar@{->>}[d] \\
&S \ar@{->>}[r] &I\ulpullback}}
\end{equation}
by taking pushout along $\pi$ and $\tau$ and pullback of the pushout diagram. Note that all the arrows are surjections except perhaps $\phi$. Note also that any pullback of surjections is also a pushout square. The proof of the following lemma is straightforward.  
 \begin{lemma} Let $\pi$ and $\tau$ be two partitions of $E$ presented as surjections as in \eqref{partsurj}.
  \begin{enumerate}[i)]
  \item The join $\pi\vee \tau$ corresponds to the surjection $E\twoheadrightarrow I$ and $\hat{0}$ is $E\twoheadrightarrow E$.
  \item The meet $\pi \wedge \tau$ corresponds to the surjection $\phi\colon E\twoheadrightarrow \im (\phi)$ and $\hat{1}$ is $E\twoheadrightarrow 1$.
 \item $\pi \wedge \tau=\hat{0}$ if and only if $\phi$ is injective.
 \item $\pi$ and $\tau$ commute if and only if $\phi$ is surjective.
 \item $\pi$ and $\tau$ are independent if and only if $\phi$ is surjective and $I=1$.
 \item For a partition $\sigma$ we have that $\pi\le \sigma$ as partitions if and only if $\sigma$ factors through $\pi$ as surjections.
 \end{enumerate}
 \end{lemma}
 Recall from \cite{Nava-Rota} that $\pi$ and $\tau$ are said to {\em commute} when, for all $p,q\in E$, if there is an $r\in E$ which belongs to the same block of $\pi$ as $p$ and to the same block of $\tau$ as $q$ then there is an $s\in E$ which belongs to the same block of $\tau$ as $p$ and to the same block of $\pi$ as $q$. Also, $\pi$ and $\tau$ are said to be {\em independent} if every block of $\pi$ meets every block of $\tau$. If two partitions are independent then they commute, and if they commute and its join is $\hat{1}$ then they are independent.
Nava and Rota \cite{Nava-Rota} define a \emph{transversal} of a given partition $\sigma$ of a set $E$ to be a pair of partitions $\pi,\tau$, such that 
\begin{enumerate}[i)]
\item $\pi\le \sigma$,
\item $\pi\wedge \tau=\hat{0}$,
\item  $\pi$ and $\tau$ commute,
\item $\pi\vee\tau=\sigma\vee \tau$.
\end{enumerate}
 In view of the lemma above, a transversal of the surjection  $\sigma\colon E\twoheadrightarrow B$ is a diagram 
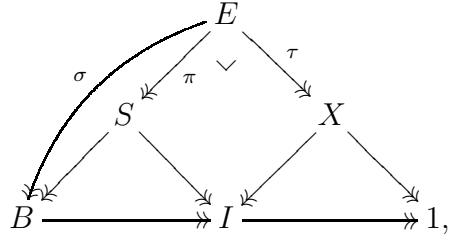
\[
\vcenter{\xymatrix{& & E \ar@{->>}[ld]^{\pi} \ar@/_1.3pc/@{->>}[lldd]_{\sigma} \ar@{->>}[rd]^{\tau}\dpullback& & \\                             
			&S\ar@{->>}[ld] \ar@{->>}[rd]& & X\ar@{->>}[ld] \ar@{->>}[rd]&\\
	     B \ar@{->>}[rr]& &I \ar@{->>}[rr]& &1,}}
\]
where the square is obtained as the pushout of $\pi$ and $\tau$. The fact that $\pi\wedge \tau=\hat{0}$ and that $\pi$ and $\tau$ commute implies that this square is also a pullback. Furthermore the condition that the pushouts $\pi\vee\tau$ and $\sigma\vee \tau$ coincide gives a map $B\twoheadrightarrow I$. Conversely, any commutative square of the form
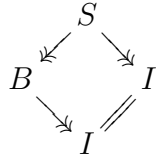
\[
\vcenter{\xymatrix@R=2ex@C=2ex{  &S\ar@{->>}[ld] \ar@{->>}[rd]&\\
	     B \ar@{->>}[rd]& &I \ar@{=}[ld]\\
	     & I&}}
\]
is a pushout in the category of surjections. Therefore the  map $B\twoheadrightarrow I$ says that $\pi\vee\tau$  coincides with $\sigma\vee \tau$.
As a consequence the groupoid of all transversals is the subgroupoid of connected objects of $T_2\sur$.

\section{Fa\`a di Bruno formula for the Green function}
\label{FdB}

In the completion of $\mathcal{P}$ define the \emph{Green function} (the terminology comes from quantum field theory \cite{GKT:FdB}) to be the series
$$A:=\sum_{\l}a_{\l}.$$
To finish, we take the opportunity to derive a formula for the comultiplication of the Green function, in close analogy with the Fa\`a di Bruno formulas of \cite{BFK:0406117} and \cite{GKT:FdB}.
\begin{proposition} Let $a_k:=\sum_{|\l|=k}a_{\l}$. Then
$$\Delta(A)=\sum_k A^k\otimes a_k.$$
\begin{proof}
This could be proved directly in $\mathcal{P}$, but we will prove it more elegantly by deriving an equivalence of groupoids related to $T\sur$ and taking homotopy cardinality. First of all, let $G$ be the inclusion $\Sg\xhookrightarrow{\;\;G\;\;} T_1\sur$ taking $a\twoheadrightarrow b$ to 
\[
\vcenter{\xymatrix@C=1ex@R=2ex{& a \ar@{->>}[ld] \ar@{->>}[rd]  & \\
             b \ar@{->>}[rr]&  &1.}}
\]
It is clear that 
$$|G|=\sum_{\l}\frac{A_{\l}}{|\aut(\l)|}=\sum_{\l}a_{\l}=A.$$
Denote by $\T$ the subgroupoid of connected objects of $T_2\sur$. Observe that $d_1^{\ast}(G)$ is precisely the inclusion $\T\xhookrightarrow{\;\;\;} T_2\sur$. Therefore  $\Delta(G)$ is the map $\T\xrightarrow{(d_2,d_0)_!d_1^{\ast}(G)} T_1\sur\times T_1\sur$.
Now, since $T\sur$ is a Segal space we have that $T_2\sur\simeq T_1\sur \times_{\B} T_1\sur$, where $\B=T_0\sur$ is the groupoid of finite sets and bijections. As a consequence $\T\simeq T_1\sur\times_{\B}\Sg$. 
 In pictures this equivalence looks like
\[
\T\simeq \left\{ \vcenter{\xymatrix{& \ast \ar@{->>}[ld] \ar@{->>}[rd]& \\
             \ast\ar@{->>}[rr]&&t_{11}, }}
 \vcenter{\xymatrix{& \ast \ar@{->>}[ld] \ar@{->>}[rd]& \\
             t_{11}\ar@{->>}[rr]&&1 }}\right\}
\]
We can decompose this as the homotopy sum of its fibres
$$\T\simeq \int^{b\in\B}(T_1\sur)_b\times{}_b\Sg,$$
where now the right $b$ subscript means $d_0$-fibre over $b$ and the left $b$ subscript means $d_1$-fibre over $b$. A fancier way to express this equation is 
$$\T\simeq \int^{b\in \B} \Sg^{b}\times \B^{b},$$
in view of the equivalences $(T_1\sur)_b\simeq \Sg^b$ and ${}_b\Sg\simeq \B^b$.

Let us now take homotopy cardinality. To compute $\Delta (A)$ we only need to know $|\Sg^b\rightarrow T_1\sur|$ and $|\B^b\rightarrow T_1\sur|$. The latter, if we call $|b|=k$, is clearly
$$|\B^b\rightarrow T_1\sur|=\sum_{|\l|=k}a_{\l}=:a_k.$$
For the former, just notice that $\Sg^b\rightarrow T_1\sur\simeq(\Sg\rightarrow T_1\sur)^{\grpdSprod k}$, since $(T_1\sur)_b\simeq (T_1\sur)_1^{+k}$. Here $\grpdSprod$ is the monoidal product in $\infgrpd_{/T_1\sur}$. This implies that
%
$$|\Sg^b\rightarrow T_1\sur|=|(\Sg\rightarrow T_1\sur)^{\grpdSprod k}|=|\Sg\rightarrow T_1\sur|^k=A^k,$$
and therefore
$$\Delta(A)=\sum_k A^k\otimes a_k,$$
as asserted. 
\end{proof}
\end{proposition}


\noindent \textsc{Departament de Matem\`atiques,\\
Universitat Aut\`onoma de Barcelona,\\
08193 Bellaterra (Barcelona),
Spain}\\
\textit{E-mail adress:} \href{mailto:acebrian@mat.uab.cat}{acebrian@mat.uab.cat},\\

\end{document}